\documentclass{amsart}
\usepackage[centertags]{amsmath}
\usepackage{amsfonts}
\usepackage{amssymb}
\usepackage{amsthm}
\usepackage{graphicx}
\usepackage{lmodern}
\usepackage{newlfont}
\usepackage[T1]{fontenc}
\usepackage{newlfont}
\usepackage{hyperref}
\usepackage{setspace}
\usepackage[english,frenchb]{babel}
\usepackage[latin1]{inputenc}
\newcommand{\ds}{\displaystyle}

\newcommand{\R}{\mathbb R}

\newcommand{\Hr}{\mathbb H}

\newtheorem{thm}{Th\'eor\`eme}[section]
\newtheorem{prop}[thm]{Proposition}
\newtheorem{lem}[thm]{Lemme}
\newtheorem{cor}[thm]{Corollaire}
\newtheorem{defn}[thm]{D\'efinition}
\newtheorem{rem}{Remarque}[section]
\numberwithin{equation}{section}

\def\md {\par \medskip}
\def\ds {\displaystyle}
\def\<{\langle}
\def\>{\rangle}

\begin{document}

\title[R\'ecurrence ou non minimalit\'e des  adh\'erences des d'orbites irr\'eguli\`eres. ]
 {R\'ecurrence ou non minimalit\'e des  adh\'erences des d'orbites irr\'eguli\`eres du flot horocyclique de finesse infinie}

\author{Masseye Gaye et Amadou Sy}
\address{Laboratoire G\'eom\'etrie et Application (LGA),
D\'epartement Math\'ematique et Informatique, UCAD-DAKAR, Senegal.}
 \email{masseye.gaye@ucad.edu.sn}
\email{amadou22.sy@ucad.edu.sn}


\begin{abstract}
La dynamique topologique du flot horocyclique $h_\mathbb{R}$ sur le fibr\'e tangent unitaire d'une surface hyperbolique  g\'eom\'etriquement finie est bien connue. En particulier sur une telle surface le flot $h_\mathbb{R}$ est minimal ou les ensembles minimaux sont les orbites p\'eriodiques. Lorsque la surface  est g\'eom\'etriquement infinie, la situation est plus complexe et  la pr\'esence d'\'eventuelles d'orbites non ferm\'ees et non denses, appel\'ees orbites irr\'eguli\`eres, rend la description des ensembles minimaux compliqu\'ee. Dans ce texte, nous allons montrer une telle orbite est r\'ecurrente ou son adh\'erence est non $h_{\R}$ minimal.  Ce qui permettrait de presque finir la description des ensembles $h_{\R}$-minimaux.
\end{abstract}
\maketitle
\selectlanguage{english}
\begin{abstract}
The topological dynamics of the horocyclic flow $h_{\R}$ on the unit tangent bundle of a geometrically finite hyperbolic surface is well known.
In particular, on such a surface, the flow $h_{\R}$ is minimal, or the minimal sets are the periodic orbits. When the surface is geometrically infinite, the situation is more complex, and the presence of possible non-closed and non-dense orbits, called irregular orbits, complicates the description of minimal sets. In this text, we will show that such an orbit is recurrent, or its closure is non-$h_{\R}$ minimal. This would allow us to almost complete the description of $h_{\R}$-minimal sets.
\end{abstract}
\selectlanguage{french}

\textbf{Mots cl\'es :} horocycli orbits, minimal set, asymptotic finness, recurrence.\\
AMS 2010 \textit{Mathematics Subject Classification.} Primary 37D40; Secondary 20H10, 14H55, 30F35.

\section{Introduction}

Dans l'\'etude des syst\`emes dynamiques, une question importante est de comprendre les ensembles minimaux, c'est-\`a-dire les ferm\'es non vide, invariants et qui ne contient aucun ferm\'e propre, non vide et invariant. Pour le flot horocyclique $h_{\R}$ sur une surface hyperbolique g\'eom\'etriquement infinie restreint \`a son ensemble non-errant $\Omega_h$, cette question est encore ouverte. Mais il est bien connu que comprendre les ensembles $h_{\R}$ minimaux revient \`a comprendre les adh\'erences des orbites irr\'eguli\`eres (non ferm\'ee et non dense dans l'ensemble non errant $\Omega_{h}$), voir par exemple \textbf{\cite{ab1}, \cite{ab2}}, \textbf{\cite{glo}, \cite{kul2}} et  \textbf{\cite{mats}}.  Pour $u\in\Omega_h$, un argument classique pour avoir la non minimalit\'e de $\overline {h_{\R}(u)}$ est de montrer que  l'ensemble $T_u=\{t\in\R^{*};\ g_t(u)\in\ \overline{h_{\R}(u)}\}\neq \emptyset ;$ c'est-\`a-dire l'orbite g\'eod\'esique $g_{\R}(u)$ rencontre l'adh\'erence de l'orbite horocyclique $h_{\R}(u)$ en un temps non nul. L'\'etude de cette intersection est li\'ee \`a la compr\'ehension de la limite inf\'erieure du rayon d'injectivit\'e le long de la demi g\'eod\'esique $u(\R_{+})$, appel\'ee finesse asymptotique de $u(\R_{+})$ et not\'ee $Inj(u(\R_{+}))$. Dans   \textbf{\cite{ab1}},\textbf{\cite{ab2}}), A. Bellis montre que :

{\it  Pour une surface hyperbolique $\Sigma$ g\'eom\'etriquement infinie, si $u$ est un \'el\'ement de $\Omega_h$  telle que $h_{\R}(u)$ non ferm\'ee et que $0\leq Inj(u(\R_{+}))<+\infty$. Alors $\overline{h_{\R}(u)}$ n'est pas $h_{\R}$-minimale.} 



Au vu du r\'esultat ci-dessus on se pose la question que passe t-il si la finesse est infinie?  

Autrement dit  la situation $Inj(u(\R_{+})) =+\infty$ donne t-elle lieu \`a de nouveaux types d'ensembles $h_{\R}$-minimaux ?  
R\'ecemment dans \textbf{\cite{dflm}}, les auteurs ont construit pour la premi\`ere fois : 

{\it  un exemple de surfaces hyperboliques g\'eom\'etriquement infinies  $\Sigma$ pour lesquelles il existe $u_0\in \Omega_h$ tel que :  

$h_{\R}(u_0)$ est  non ferm\'ee , non dense, de finesse asymptotique $Inj(u_0(\R_{+}))=+\infty$ et  $\overline{h_{\R}(u)}$  est minimale.}

Dans ce papier, nous montrons le r\'esultat suivant qui est une  dichotomie entre la r\'eccurence des orbites irr\'eguli\`eres de finesse infinie et la non minimalit\'e de leur adh\'erence.
\begin{thm}\label{pa}
Soit $\Sigma=\Gamma\setminus\mathbb{H}$ une surface hyperbolique  g\'eom\'etriquement infinie  o\'u $\Gamma$ est un groupe fuschien sans \'el\'ement elliptique. 
Soit $u_0\in \Omega_h$ tel que $h_{\R}(u_0)$ est  irr\'eguli\`ere (non ferm\'e et non dense) et de finesse asymptotique $Inj(u_0(\R_{+}))=+\infty$.
Alors :

- soit $h_{\R}(u_0)$ est r\'ecurrente ;

- soit l'adh\'erence $\overline{h_{\R}(u_0)}$ n'est pas $h_{\R}$ -minimale.
\end{thm}

Puisqu'un ensemble $h_{\R}$-minimal est l'adh\'erence d'une orbite, les r\'esultats du  th\'eor\`eme \ref{pa} combin\'es avec ceux de Bellis dans \textbf{\cite{ab2}}), nous dit les \'eventuels ensemble minimaux pour le flot horocyclique sur une surface hyperbolique g\'eom\'etriquement infinie restreint \`a son ensemble non-errent sont ses orbites ferm\'ees ou l'adh\'erence de ses orbites irr\'eguli\'eres non r\'ecurrentes de finesse infinie.

 Question : {\it l'adh\'erence d'une orbite irr\'eguli\`ere r\'ecurrente de finesse infine est -elle toujours minimale}?



La suite du papier est organis\'ee comme suit : la premi\`ere section est resev\'ee aux pr\'eliminaires dont nous aurons besoin dans la deuxi\`eme section consacr\'es \`a la d\'emonstration du th\'eor\`eme \ref{pa}. 

\section{Pr\'eliminaires}

Dans cette section, nous rappelons seulement les notions et les r\'esultats dont nous aurons besoin pour la d\'emonstration du th\'eor\`eme \ref{pa}.

\subsection{Nature des points de l'ensemble limite.}
Soit $\Sigma=\Gamma\setminus\mathbb{H}$ un surface hyperbolique  o\'u $\Gamma$ est un groupe fuschien sans \'el\'ement elliptique. Notons $\Lambda=\overline{\Gamma.i}\setminus\Gamma.i$ l'ensemble limite du groupe $\Gamma.$ Dans  $\Lambda$ nous distinguons quatre types de points selon l'orbite $\Gamma.i$ intersetcte les horidisques ouverts bas\'e en ce point limite (voir \textbf{\cite{dalstar}}).
\begin{itemize}
  \item un point $\xi$ de $\Lambda$ est dit \textit{horocyclique} si pour tout horodisque ouvert $\mathcal{O}_\xi$ bas\'e en $\xi$, l'ensemble $\Gamma.i\cap \mathcal{O}_\xi$ est infini. On note $\Lambda_h$ l'ensemble de ces points.
  \item un point $\xi$ de $\Lambda$ est dit \textit{discret} si pour tout horodisque ouvert $\mathcal{O}_\xi$ bas\'es, l'ensemble $\Gamma.i\cap \mathcal{O}_\xi$ est fini. On note $\Lambda_d$ l'ensemble de ces points.
  \item un point $\xi$ de $\Lambda$ est dit \textit{parabolique} s'il est fix\'e par une isom\'etrie parabolique de $\Gamma$. On note $\Lambda_p$ l'ensemble de ces points.
  \item un point $\xi$ de $\Lambda$ est dit \textit{irr\'egulier} s'il n'est ni horocyclique, ni discret et ni parabolique. On note $\Lambda_{irr}$ l'ensemble de ces points.
\end{itemize}
Les ensembles $\Lambda_h$, $\Lambda_d$, $\Lambda_p$ et $\Lambda_{irr}$ forment une partition $\Gamma$ invariante de $\Lambda.$
\begin{rem} \label{plhd} Soit  $\infty$ est un point de l'ensemble limite de $\Lambda$.

- $\infty \notin\Delta_h $ si et seulement s'il existe $M>0$ tel pour tout $\gamma\in\Gamma, Im(\gamma(i))\leq M$.

- $\infty\notin\Delta_d$ s'il existe une suite non constante $(Im(\gamma_n(i)))_{n\geq 1}$ qui tend vers $l>0.$
\end{rem}
Le th\'eor\`me suivant bien connu (voir par exemple \textbf{\cite{dal1}})  donne une caract\'erisation de la finitude de la surface en fonction des points limites
\begin{thm}
La surface $\Sigma=\Gamma\setminus\Hr$ est g\'eom\'etriquement finie si et seulement si $\Lambda=\lambda_h\cup\Lambda_p.$
\end{thm}
\subsection{Birapport et intersection de g\'eod\'esiques}
La r\'ef\'erence que nous utilisons pour cette sous section est \cite{fl}.
\begin{defn}
Soient $a, b, c \;et\; d $ quatre points distincts du bord de $\mathbb{H}$. Le birapport de ces points est d\'efini par :
$$
[a; b; c; d] = \frac{(a - c)(b - d)}{(a - d)(b - c)} \;avec \;la\; convention\; [a; b; c; \infty] = \frac{a - c}{b - c}.
$$
\end{defn}
\begin{lem}{{Voir\cite{fl}; prop $2.1.16$}}
Soient $(a,b)$ et $(c,d)$ deux g\'eod\'esiques s'intersectant en un point $x$ . Si les points du bord  sont ordonn\'es comme  suit $(a; c; b; d)$, alors l'angle  $\beta \in [0,\pi]$ entre les g\'eod\'esiques $(a,b)$ et $(c,d)$ v\'erifie :
$$
[a; c; d; b] = \frac{\cos(\beta) + 1}{2}.
$$
\end{lem}

\subsection{Nature topologique des orbites du flot horocyclique.} Sur une surface hyperbolique g\'eom\'etriquement finie, la dynamique topologique du flot horocyclique $h_{\R}$
 est bien connue : {\it les orbites sont denses ou p\'eriodique. En particulier le flot $h_\mathbb{R}$ est minimal ou les ensembles $h_{\R}$-minimaux sont ses orbites p\'eriodiques.}

Lorsque la surface hyperbolique est g\'eom\'etriquement infinie, la dynamique topologique  est plus compliqu\'ee et d\'epend de l'action de $\Gamma$ sur le bord \`a l'infini $\partial\mathbb{H}$. L'ensemble des orbites de $h_\mathbb{R}$ se divise en quatre sous ensemble de nature topologique diff\'erente :
\begin{itemize}
\item l'ensemble des orbites denses dans $\Omega_h$;

\item l'ensemble des orbites p\'eriodique (compactes);

\item l'ensemble des orbites ferm\'es non compactes;

\item l'ensemble des orbites irr\'eguli\`eres (ni ferm\'ee, ni dense dans $\Omega_h$);
\end{itemize}
et il existe toujours une orbites qui n'est ni dense et ni p\'eriodique, voir \textbf{\cite{star1}}, \textbf{\cite{dalstar}}. Il est donc naturel de se demander que devient la rigidit\'e du flot $h_\mathbb{R}$ observ\'ee dans le cadre g\'eom\'etriquement finie. Nous remarquons alors :
\begin{itemize}
  \item {\it le flot $h_{\R}$ n'est pas minimal, c'est-\`a-dire $\Omega_h$ n'est pas un ensemble minimal;}
  \item {\it les orbites ferm\'ees de $h_{\R}$ sont des ensembles minimaux pour $h_{\R}$;} 
  \item {\it les \'eventuels ensembles $h_{\R}$-minimaux diff\'erent orbites ferm\'ees, sont les adh\'erences de ses orbites irr\'eguli\`eres.}  
  \end{itemize}
Jusqu'\`a pr\'esent toutes les investigations men\'ees sur une large classe de surfaces hyperboliques g\'eom\'etrqiement infinies, montrent que les adh\'erences des orbites irr\'eguli\`eres ne sont pas des ensembles $h_{\R}$-minimaux (voir par exemple \textbf{\cite{glo}},  \textbf{\cite{kul2}}, \textbf{\cite{mats}}). 

\subsection{Correspondance entre point limite et nature topologique des orbites.}
Pour $u\in T^1\Sigma,$ nous notons $\widetilde u$ un relev\'e de $u$ et $u^{+}$ l'extr\'emit\'e positive de la g\'eod\'esique de $\Hr$ d\'efinie par $\widetilde u.$ Il est bien connu que $u$ appartient \`a $\Omega_h$ si et seulement si l'extr\'emit\'e $u^+$ appartient \`a l'ensemble $\Lambda$. Il existe donc une correspondance entre la nature topologique de l'orbite $h_\mathbb{R}(u)$ et la nature du point limite $u^+$ appartenant \`a $\Lambda$ (voir par exemple \textbf{\cite{star1}} et \textbf{\cite{dalstar}}). Plus pr\'ecis\'ement nous avons la correspondance suivante :
{\it \begin{itemize}
  \item l'orbite $h_\mathbb{R}(u)$ est dense dans $\Omega_h$ si et seulement si $u^+\in\Lambda_h$,
  \item l'orbite $h_\mathbb{R}(u)$ est p\'eriodique si et seulement si $u^+\in\Lambda_p$,
  \item l'orbite $h_\mathbb{R}(u)$ est ferm\'ee et non p\'eriodique si et seulement si $u^+\in\Lambda_d$,
  \item l'orbite $h_\mathbb{R}(u)$ est irr\'eguli\`ere (orbite qui n'est ni ferm\'ee, ni dense dans $\Omega_h$) si et seulement si $u^+\in\Lambda_{irr}$.
\end{itemize}}

\subsection{Convergence dans $T^1\Sigma$}
Il s'agit ici de rappeler des r\'esultats sur la convergence dans $T^1\Sigma$ qui seront n\'ecessaire pour d\'emontrer le th\'eor\`eme \ref{pa}. Pour $x\in \Sigma,\ u\in\ T^1\Sigma$ notons respectivement $\widetilde x$ et $\widetilde u$ des relev\'es quelconques de $x$ et de $u$ dans $\Hr$ et dans $T^1\Hr.$ Notons \'egalement par $(\widetilde u(t))_{t\geq 0}$ la param\'etrisation par unit\'e de loungueur du demi rayon g\'eod\'esique $\widetilde u(\R_{+})$ de $\Hr$, d\'efini par $\widetilde u$ d'origine $\widetilde u(0)$ et d'extr\'emit\'e $\widetilde u(\infty).$ D\'esignons aussi par $u(\R_{+})$ le d\'emi- rayon g\'eod\'esique de $\Sigma,$ d\'efinie par $u$, param\'etr\'e par $(u(t))_{t\geq 0}$ le projet\'e de $(\widetilde u(t))_{t\geq 0}$ dans $\Sigma$. Pour $z$ dans $\Hr$, notons $B_{\widetilde u(\infty)}(\widetilde u(0), z)$ le cocycle de Busemann centr\'e en $\widetilde u(\infty)$, calcul\'e en $\widetilde u(0), z$.
L'ensemble $\{ z\in\Hr, B_{\widetilde u(\infty)}(\widetilde u(0), z)=0\}$ est l'horocycle centre en $\widetilde u(\infty)$ et passant par $\widetilde u(0)$  (pour les propri\'et\'es sur le cocylce de Busemann on peut voir \textbf{\cite{dal1}}).
\`A pr\'esent nous allons voir que montrer qu'un r\'eel $t$ est dans l'en semble $T_u$ , est \'equivalent \`a construire une suite $(\gamma_n)_{n\geq 0}$ d'\'el\'ements de $\Gamma$, qui v\'erifie deux propri\'et\'es. Cela repose sur la convergence horocyclique d\'emontr\'ee dans \textbf{\cite{ab1}}). Plus pr\'ecis\'ement nous avons la proposition suivante :
\begin{prop} (Dans \textbf{\cite{ab1}}, Chap $2$, Prop $2.3.2$)
  Soient $u$ et $v$ deux \'el\'ements $T^1\Sigma$. Alors $v$ appartient \`a l'adh\'erence  $\overline{h_{\R}(u)}$ si et seulement s'il existe une suite $(\gamma_n)_{n\geq 0}$ d'\'el\'ements de $\Gamma$ telle que les deux conditions suivantes soient v\'erifi\'ees :
\md
  $(i)\ \ds\lim_{n\rightarrow +\infty}\gamma_n\widetilde u(\infty)=\widetilde v(\infty).$ 
\md

  $(ii)\ \ds\lim_{n\rightarrow +\infty}B_{\widetilde u(\infty)}(\gamma_n^{-1}i,\widetilde u(0))=B_{\widetilde v(\infty)}(i,\widetilde v(0)).$
\end{prop}

Pour $t\neq 0,$ en posant $v=g_tu,$ alors il existe $\alpha\in\Gamma$ tel que $\widetilde v(\infty)=\alpha\widetilde u(\infty)$ et que $\widetilde v(0))=\widetilde g_t\alpha\widetilde u(0).$  Les points $(i)$ et $(ii)$ de la proposition pr\'ec\'edente deviennent alors :
$$(i)'\ds\lim_{n\rightarrow +\infty}\gamma_n\widetilde u(\infty)=\alpha\widetilde u(\infty).$$ 
$$(ii)'\ \ds\lim_{n\rightarrow +\infty}B_{\widetilde u(\infty)}(\gamma_n^{-1}i,\widetilde u(0))=B_{\alpha\widetilde u(\infty)}(i,\alpha\widetilde g_t\widetilde u(0))=B_{\widetilde u(\infty)}(\alpha^{-1}i,\widetilde u(0))+t.$$ Nous obtenons alors :
\begin{cor} \label{tunv} Soit $u$ un \'el\'ement de $T^1\Sigma$. Alors un r\'eel $t$ non nul est dans l'ensemble $T_u$ si et seulement s'ils existent $\alpha\in\ \Gamma$ et une suite $(\gamma_n)_{n\geq 0}$ d'\'el\'ements de $\Gamma$ telle que les deux conditions suivantes soient v\'erifi\'ees :
\md
  $(i)\ \ds\lim_{n\rightarrow +\infty}\gamma_n\widetilde u(\infty)=\alpha\widetilde u(\infty).$ 
\md

  $(ii)\ \ds\lim_{n\rightarrow +\infty}B_{\widetilde u(\infty)}(\gamma_n^{-1}i,\alpha^{-1}i)=t$

\end{cor}

\begin{cor} \label{or} Soit $u$ un \'el\'ement de $T^1\Sigma$. Alors l'orbite $h_{\R}(u)$ est r\'ecurrente si et seulement s'il existe tune suite $(\gamma_n)_{n\geq 0}$ d'\'el\'ements de $\Gamma$ tous distincts telle que les deux conditions suivantes soient v\'erifi\'ees :
\md
  $(i)\ \ds\lim_{n\rightarrow +\infty}\gamma_n\widetilde u(\infty)=\widetilde u(\infty).$ 
\md

  $(ii)\ \ds\lim_{n\rightarrow +\infty}B_{\widetilde u(\infty)}(\gamma_n^{-1}i,i)=0$

\end{cor}
\section{D\'emonstration du th\'eor\`ereme \ref{pa}}
Soit $\Sigma=\Gamma\setminus\mathbb{H}$ une surface hyperbolique  g\'eom\'etriquement infinie  o\'u $\Gamma$ est un groupe fuschien sans \'el\'ement elliptique. 
Soit $u_0\in \Omega_h$ tel que $h_{\R}(u_0)$ est  irr\'eguli\`ere (non ferm\'e et non dense) et de finesse asymptotique $Inj(u_0(\R_{+}))=+\infty$. La d\'emonstration du th\'eor\`eme \ref{pa} se fait en plusieurs  \'etapes. Elle consiste \`a construire une suite $(\gamma_n)_{n\geq 0}$ d'\'el\'ements de $\Gamma$ qui permet d'appliquer le corollaire \ref{tunv} ou le corollaire \ref{or}. Placons nous dans le demi-plan de Poincar\'e $\mathbb{H}$ et quitte \'a conjuguer $\Gamma$ par un \'el\'ement de $PSL_2(\R)$, nous suposerons $u_0(\R_+)$ est le projet\'e sur $\Sigma$ de la demi-g\'eodisique orient\'ee $[i\infty)$ de $\mathbb{H}$. Ainsi dans la suite de ce texte, $u_0$ le projet\'e sur de l \'el\'ement $\widetilde u_0$ de $T^1\mathbb H$ bas\'e en $i$ et pointant vers $\infty$ qui est un point limite irr\'egulier du groupe $\Gamma.$
\begin{lem} \label{lem1}Soit $\Gamma$ un groupe fuschien g\'eom\'etriquement infinie sans \'el\'ement elliptique et $\infty\in\Lambda_{irr}$. Alors ils existent deux r\'eels $M>m>0$ et une suite $(\gamma_n)_{n\geq 1}$ d'\'el\'ements de $\Gamma$ telle que la suite $(\gamma_n(i))_{n\geq 1}$ est non constante et v\'erifie :

$$m\leq Im(\gamma_n(i))\leq M,\ \ds\lim_{n\rightarrow +\infty}\gamma_n(i)=\infty\ \mbox{et}\ \ \ds\lim_{n\rightarrow +\infty}\gamma_n(\infty)=\infty$$

\end{lem}
{\it \bf {D\'emonstration} :} Puisque $\infty$ n'est pas un point limite horocyclique ni un point limite discret, pour le groupe $\Gamma,$ d'apr\`es la remarque \ref{plhd}, ils existent $M>m>0$ et une suite $(\alpha_n)_{n\geq 1}$ d'\'el\'ements de $\Gamma$ telle que la suite $(Im(\alpha_n(i)))_{n\geq 1}$ est non constante et v\'erifie $m\leq Im(\alpha_n(i))\leq M$.
 Comme la suite $(\alpha_n(i))_{n\geq 1}$ n'est pas constante \'a partir d'un certain rang et que $0<m\leq Im(\alpha_n(i))$, alors $\infty$ est l'unique  point d'accumalation de la suite
 $(\alpha_n(i))_{n\geq 1}$. Par suite il existe un sous suite $(\gamma_n(i))_{n\geq 1}$ de la suite $(\alpha_n(i))_{n\geq 1}$ qui v\'erifie : $$m\leq Im(\gamma_n(i))\leq M \ \mbox{et}\ \ds\lim_{n\rightarrow +\infty}\gamma_n(i)=\infty\ .$$ Montrons \`a pr\'esent que $\gamma_n(\infty)$ tend vers $\infty$. Pour tout entier $n\geq 1,$ pour tout $z$ dans $\mathbb H,$ posons: $\gamma_n(z)=\ds\frac{a_nz+b_n}{c_nz+d_n}$ avec $a_n, b_n, c_n$ et $d_n$ dans $\R$ et  $a_nd_n-b_nc_n=1$. Comme le point $\infty$ n'est fix\'e par aucun \'el\'ement de $\Gamma\setminus\{Id\}$, alors pour tout entier $n\geq 1, c_n\neq 0$ et  $\gamma_n(\infty)=\ds\frac{a_n}{c_n}$. Pour tout entier $n\geq 1,$ nous avons  \'egalement $m\leq Im(\gamma_n(i))=\ds\frac{1}{c_n^2+d_n^2}\leq M.$ Il en r\'esulte que lees suites $(c_n)_{n\geq 1}$ et $(d_n)_{n\geq 1}$ sont born\'ees et quitte \`a extraire sumultan\'ement deux sous suites not\'ees encore $(c_n)_{n\geq 1}$ et $(d_n)_{n\geq 1}$ , alors elles convergent  respectivement vers $c$ et $d$ v\'erifiant $c^2+d^2\geq\ds\frac{1}{M}>0.$  
 
 Par un calcul direct nous avons \'egalement $\Re(\gamma_n(i))=\ds\frac{a_n}{c_n}-\ds\frac{d_n}{c_n(d_n^2+c_n^2)}$. Si $c\neq 0,$ la suite $\ds\frac{d_n}{c_n(d_n^2+c_n^2)}$ converge vers $\ds\frac{d}{c(d^2+c^2)}$. Comme $\gamma_n(i)$ tend vers $\infty$ et $Im(\gamma_n(i))$ est born\'es alors $\Re(\gamma_n(i))$ tend vers $\infty$; par suite $\gamma_n(\infty)=\ds\frac{a_n}{c_n}$ tend vers $\infty$. Si $c=0,$ quitte \`a consid\'erer une sous suite, nous supposons que la suite $(a_n)_{n\geq 1}$ admet une limte $a\in\R^{*}\cup\{-\infty,+\infty\}$ car si $(a_n)_{n\geq 1}$ convergeait vers $0,$ nous aurions $Im(\gamma_n^{-1}(i))=\ds\frac{1}{c_n^2+a_n^2}$ qui tenderait vers $+\infty$ et nous serions en contradiction avec $\infty$ est un point limite discret. Nous en d\'eduisons $\gamma_n(\infty)$ tend vers $\infty$.
  $$\hfil\Box$$

Dans toute la suite, nous d\'esignerons par $(\gamma_n)_{n\geq 1}$  la suite donn\'ee dans le lemme \ref{lem1} et par  $a_n, b_n, c_n$ et $d_n$ les coefficients associ\'ees \`a l'isom\'etrie $\gamma_n$. Alors les suites $(c_n)_{\geq 1}, (d_n)_{\geq 1}$ convergent respectivement vers $c$ et $d$ v\'erifiant $c^2+d^2>0$ et la suite $(a_n)_{n\geq 1}$ admet une limite $a\in \R\cup\{-\infty\,+\infty\}$.

\begin{lem} \label{lem2}Soient $\Gamma$ un groupe fuschien g\'eom\'etriquement infinie sans \'el\'ement elliptique et $u_0\in \Omega_h$ tel que $h_{\R}(u_0)$ est  irr\'eguli\`ere (non ferm\'e et non dense) et de finesse asymptotique $Inj(u_0(\R_{+}))=+\infty$. Alors $(a_n)_{n\geq 1}$ admet $a\in\{-\infty,+\infty\}$ pour limite lorsque $n$ tend vers $+\infty$.

\end{lem}
{\it \bf {D\'emonstration} :} Supposons que   $(a_n)_{n\geq 1}$ admet $a\in\R$ pour limite lorsque $n$ tend vers $+\infty$. Comme les suites  $(c_n)_{\geq 1}, (d_n)_{\geq 1}$ convergent respectivement vers $c$ et $d$ et $b_nc_n=a_nd_n-1$ alors la suite $(b_nc_n)_{n\geq 1}$ converge vers $ad-1.$ Par ailleurs nous avons $ |\gamma_n(i)|^2=\ds\frac{a_n^2+b_n^2}{c_n^2+d_n^2}$ qui tend vers $+\infty$ d'apr\`es le lemme \ref{lem1} et comme les suites $(a_n)_{n\geq 1}, (c_n)_{\geq 1}, (d_n)_{\geq 1}$ convergent respectivement vers $a, c, d$ avec $c^2+d^2>0$;  nous en d\'eduisons que la suite $(b_n^2)_{\geq 1}$ tend vers $+\infty.$

D'une part, comme par hypoth\`ese  $Inj(u_0(\R_{+}))=+\infty$, en notant $(ie^t)_{t\geq 0}$ le param\'etrage par unit\'e de longueur de la demi g\'eod\'esique $\widetilde u_0=[i\infty)$ de $\mathbb H$, pour toute suite de temps $(t_n)_{n\geq}$ tendant vers $+\infty$, nous avons : $$\ds\lim_{n\rightarrow +\infty}\inf_{\gamma\in\Gamma\setminus\{Id\}}d(ie^{t_n},\gamma(ie^{t_n}))=+\infty;\ \mbox{en particulier}\ \lim_{n\rightarrow +\infty}d(ie^{t_n},\gamma_n(ie^{t_n}))=+\infty.$$
D'autre part nous avons $d(ie^{t_n},\gamma_n(ie^{t_n}))=2 Argsh\bigg(\ds\frac{|ie^{t_n}-\gamma_n(ie^{t_n})|} {2\sqrt{e^{t_n}Im(\gamma_n(ie^{t_n}))}}\bigg),$ ce qui donne :
$$ d(ie^{t_n},\gamma_n(ie^{t_n}))=2 Argsh\bigg(\ds\frac{\sqrt{b_n^2 e^{-2t_n}+c_n^2 e^{2t_n}+d_n^2+a_n^2-2}}{2}\bigg).$$
Pour $n\geq 1,$ en posant $t_n=\ln(b_n^2),$ nous obtenons :
$$ d(ie^{t_n},\gamma_n(ie^{t_n}))=2 Argsh\bigg(\ds\frac{\sqrt{b_n^2 c_n^2 +d_n^2+a_n^2-1}}{2}\bigg).$$ Par passage \`a la limite lorsque lorsque $n$ tend vers $+\infty,$ nous avons : $$\ds\lim_{n\rightarrow +\infty}d(ie^{t_n},\gamma_n(ie^{t_n}))=2 Argsh\bigg(\ds\frac{\sqrt{(ad-1)^2 +d^2+a^2-1}}{2}\bigg);$$ ce qui est en contradiction  avec  $\ds\lim_{n\rightarrow +\infty}d(ie^{t_n},\gamma_n(ie^{t_n}))=+\infty.$
$$\hfil\Box$$

\begin{prop} \label{prop1}Soient $\Gamma$ un groupe fuschien g\'eom\'etriquement infinie sans \'el\'ement elliptique et $u_0\in \Omega_h$ tel que $h_{\R}(u_0)$ est  irr\'eguli\`ere (non ferm\'e et non dense) et de finesse asymptotique $Inj(u_0(\R_{+}))=+\infty$. Alors $(c_n)_{n\geq 1}$ converge vers $0$.

\end{prop}
{\it \bf {D\'emonstration} :}  D'apr\`es les lemmes \ref{lem1} et \ref{lem2}, les suites $(a_n){n\geq 1}, (c_n){n\geq 1}$ et $(d_n){n\geq 1}$ admettent pour limites respectives $a\in\{-\infty,+\infty\}, c$ et $d$, avec $c^2+d^2>0.$ Montrons par absurde que $c=0$ et pour cela supposons $c\neq 0$. Pour tout $n\geq 1,$ notons  $X_n=\ds\frac{a_n-d_n+\sqrt{(a_n+d_n)^2-4}}{2c_n}$ et $Y_n=\ds\frac{a_n-d_n-\sqrt{(a_n+d_n)^2-4}}{2c_n}$ les points fixes de l'isom\'etrie $\gamma_n$. Puisque $(a_n){n\geq 1}$ admet $\pm\infty$ pour limite et $(c_n){n\geq 1}$ converge vers $c\neq 0,$ alors $X_n\sim\ds\frac{a_n}{c_n}$ et donc $(X_n){n\geq 1}$ admet $+ \infty$ ou $-\infty$ pour limte.  En remarquant  $Y_n=\ds\frac{2-2a_nd_n}{c_n(a_n-d_n+\sqrt{(a_n+d_n)^2-4})}$, nous avons $(Y_n)_{n\geq 1}$ qui converge vers $\ds\frac{-d}{c}.$  

Supposons maintenant $(X_n){n\geq 1}$ tend vers $+\infty$. Le cas $(X_n){n\geq 1}$ tend vers $-\infty$ se traite de la m\^eme mani\`ere. Soient $t\in]0,\ds\frac{1}{2}[$ fix\'e et $Z_n=tX_n$. Pour tout entier $n\geq 1$ assez grand, nous avons $0<Z_n<X_n$ et $Y_n<Z_n<X_n$. Notons $(\beta_nZ_n)$ l'unique g\'eod\'esique de $\mathbb H$ orthogonal \`a la g\'eod\'esique $(Y_nX_n).$ Nous avons alors $[Y_n, X_n, \beta_n, Z_n]=-1$ et donc $\beta_n=\ds\frac{X_n(tX_n+tY_n-2Y_n)}{(2t-1)X_n-Y_n}\sim\ds\frac{tX_n}{2t-1}<0$ et tend vers $-\infty$. Ainsi pour $n$ assez grand, nous avons $\beta_n<0<Z_n$ et par cons\'equent les g\'eod\'esiques orient\'ees $(\beta_n Z_n)$ et $(0\infty)$ s'intersectent dans $\mathbb H.$ En notant $\theta_n$ l'angle entre les g\'eod\'esiques $(\beta_n Z_n)$ et $(0\infty)$, nous obtenons : $$\cos\theta_n=2[\beta_n,0,\infty,Z_n]-1=\frac{2Z_n}{Z_n-\beta_n}-1=\frac{\beta_n+Z_n}{Z_n-\beta_n}= \frac{1+\frac{Z_n}{\beta_n}}{-1+\frac{Z_n}{\beta_n}}.$$ Comme $\ds\lim_{n\rightarrow +\infty}\ds\frac{Z_n}{\beta_n}=2t-1,$ nous obtenons alors $\ds\lim_{n\rightarrow +\infty}\cos\theta_n=\ds\frac{t}{t-1}$.

{\it Supppons  $-\ds\frac{d}{c}<0$.} 

\begin{figure}
    \centering
    \includegraphics[width=0.9\linewidth]{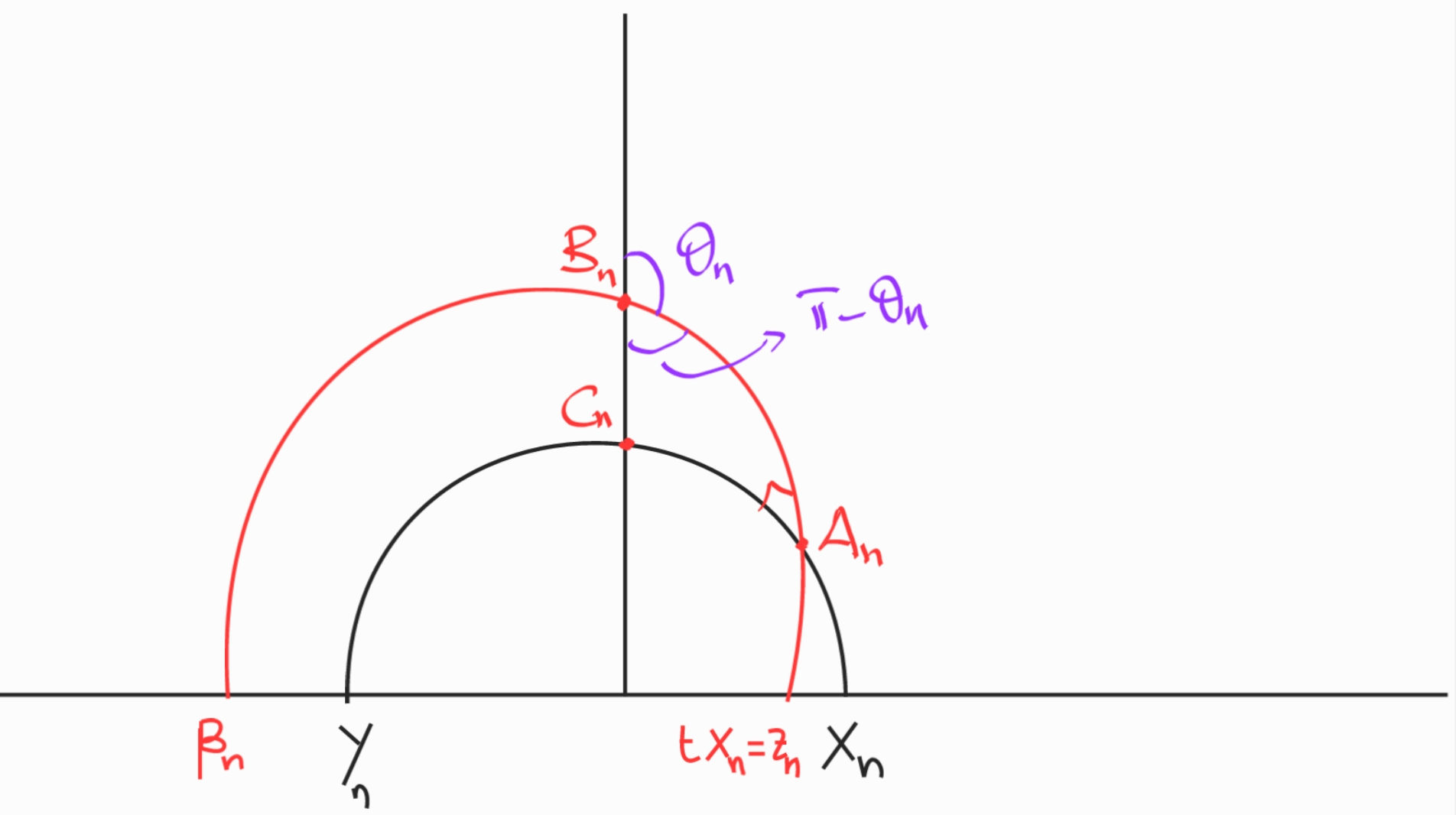}
    \caption{Triangle Hyperbolique $A_{n}B_{n}C_{n}$}
    \label{TRi}
\end{figure}

Comme $(Y_n)_{n\geq 1}$ converge vers  $-\ds\frac{d}{c}$, alors pour $n$ assez grand $Y_n<0$ et donc $\beta_n<Y_n<0<Z_n<X_n<\infty$. Par suite les g\'eod\'esiques $(Y_nX_n)$ et $(0\infty)$ s'intersectent en $C_n=i\sqrt{-Y_nX_n}$ et  les g\'eod\'esiques $(\beta_nZ_n)$ et $(0\infty)$ s'intersectent en $B_n=i\sqrt{-\beta_nZ_n}$. En notant $A_n=u_n+iv_n$ le point d'intersection des g\'eod\'esique $(Y_nX_n)$ et $(\beta_nZ_n)$, le triangle hyperbolique $A_nB_nC_n$ est rectangle en $A_n$ et l'angle au sommet $B_n$ vaut $\pi-\theta_n$ (voir figure \ref{TRi})
. En notant $l_n$ la distance hyperbolique entre $A_n$ et $B_n$ et $k_n$ celle entre $B_n$ et $C_n$, d'apr\`es les formules trigonom\'etriques hyperboliques sur un triangle hyperbolique rectangle, nous avons :
$$\cos(\pi-\theta_n)=-\cos\theta_n=\ds\frac{\tanh l_n}{\tanh k_n}.$$

Mais comme $k_n=\ln\bigg(\sqrt{\ds\frac{Z_n\beta_n}{Y_nX_n}}\bigg)=\ln\bigg(\sqrt{\ds\frac{t\beta_n}{Y_n}}\bigg)$ alors $\ds\lim_{n\rightarrow+\infty}k_n=+\infty$ et par suite $\ds\lim_{n\rightarrow+\infty}\tanh k_n=1$.

Nous avons aussi : $l_n=d(A_n,B_n)=2 Argsh\bigg(\ds\frac{|A_n-B_n|} {2\sqrt{v_n\sqrt{-\beta_nZ_n}}}\bigg)$. \\Or :

$$\frac{|A_n-B_n|^2} {4v_n\sqrt{-\beta_nZ_n}}=\frac{u_n^2+\bigg(v_n-\sqrt{-\beta_nZ_n}\bigg)^2}{4v_n\sqrt{-\beta_nZ_n}}=\frac{u_n^2+v_n^2-2v_n\sqrt{-\beta_nZ_n}-\beta_nZ_n}{4v_n\sqrt{-\beta_nZ_n}}.$$

Comme $A_n$ appartient \`a la g\'eod\'esique $(\beta_nZ_n)$ qui n'est rien d'autre que le cercle euclidien centr\'e sur l'axe r\'eel en $\ds\frac{\beta_n+Z_n}{2}$ et de rayon $\ds\frac{Z_n-\beta_n}{2}$, alors nous obtenons $u_n^2+v_n^2=u_n(Z_n+\beta_n)-\beta_nZ_n$. Par suite nous  d\'eduisons :

$$\frac{|A_n-B_n|^2} {4v_n\sqrt{-\beta_nZ_n}}=\frac{u_n}{4v_n}\bigg(\sqrt{\frac{-Z_n}{\beta_n}}- \sqrt{\frac{-\beta_n}{Z_n}}\bigg)+\frac{\sqrt{-\beta_nZ_n}}{2v_n}-\frac{1}{2}.$$

En traduisant $A_n$ est un point d'intersection des g\'eod\'esiques $(Y_nX_n)$ et $(\beta_nZ_n)$, nous obtenons :
$$u_n=\frac{\beta_nZ_n-X_nY_n}{\beta_n+Z_n-X_n-Y_n}= \frac{t\beta_nX_n-X_nY_n}{\beta_n+tX_n-X_n-Y_n}=\frac{t-\frac{Y_n}{\beta_n}}{\frac{1}{X_n}+\frac{t-1}{\beta_n}-\frac{Y_n}{X_n\beta_n}} .$$
Comme $0<t<\ds\frac{1}{2},\ Z_n=tX_n, \ds\lim_{n\rightarrow+\infty}\beta_n=-\infty$ et $ \ds\lim_{n\rightarrow+\infty}\frac{\beta_n}{X_n}=\ds\frac{t}{2t-1}$, nous avons les limites suivantes :

$$\lim_{n\rightarrow+\infty}u_n=+\infty, \lim_{n\rightarrow+\infty}\frac{u_n}{X_n}=\frac{t^2}{(1-t)^2},  \lim_{n\rightarrow+\infty}\frac{u_n}{Z_n}=\frac{t}{(1-t)^2}, \lim_{n\rightarrow+\infty}\frac{u_n}{\beta_n}=\frac{t(2t-1)}{(1-t)^2}.$$
En traduisant $A_n$ appartient \`a la g\'eod\'esique $(Y_nX_n)$, par un calcul direct nous avons $v_n^2=-u_n^2+u_n(X_n+Y_n)-X_nY_n$ et donc :

$$\lim_{n\rightarrow+\infty}\frac{v_n^2}{u_n^2}=\lim_{n\rightarrow+\infty}-1+\frac{X_n+Y_n}{u_n}-\frac{X_nY_n}{u_n^2}=\frac{1-2t}{t^2}\ \mbox{et}\ \lim_{n\rightarrow+\infty}\frac{u_n}{v_n}=\frac{t}{\sqrt{1-2t}};$$ d'o\`u $\ds\lim_{n\rightarrow+\infty}\frac{u_n}{2v_n}\bigg(\sqrt{\frac{-Z_n}{\beta_n}}- \sqrt{\frac{-\beta_n}{Z_n}}\bigg)=\ds\frac{-t^2}{1-2t}$. Mais nous avons aussi 

$$\frac{v_n^2}{-\beta_nZ_n}=\frac{u_n^2}{\beta_nZ_n}-\frac{u_n(X_n+Y_n)}{\beta_nZ_n}+\frac{X_nY_n}{\beta_nZ_n}= \frac{u_n^2}{\beta_nZ_n}-\frac{u_n(X_n+Y_n)}{t\beta_nX_n}+\frac{Y_n}{t\beta_n},$$  d'o\`u :
$$\ds\lim_{n\rightarrow+\infty}\frac{v_n^2}{-\beta_nZ_n}=\frac{(1-2t)^2}{(t-1)^4}\ \mbox{et}\ \lim_{n\rightarrow+\infty}\frac{\sqrt{-\beta_nZ_n}}{v_n}=\frac{(t-1)^2}{1-2t}.$$
Finalement nous obtenons  $\ds\lim_{n\rightarrow+\infty}\frac{|A_n-B_n|^2} {4v_n\sqrt{-\beta_nZ_n}}=0$ et donc $\ds\lim_{n\rightarrow+\infty}l_n=0$.
Ainsi $\ds\lim_{n\rightarrow+\infty}\cos\theta_n=\ds\lim_{n\rightarrow+\infty}\ds\frac{-\tanh l_n}{\tanh k_n}=0\neq\ds\frac{t}{t-1}=\ds\lim_{n\rightarrow+\infty}\cos\theta_n$  car $0<t<\ds\frac{1}{2}$
et  ceci est absurde.

{\it Supposons $-\ds\frac{d}{c}>0$}

\begin{figure}
    \centering
    \includegraphics[width=0.9\linewidth]{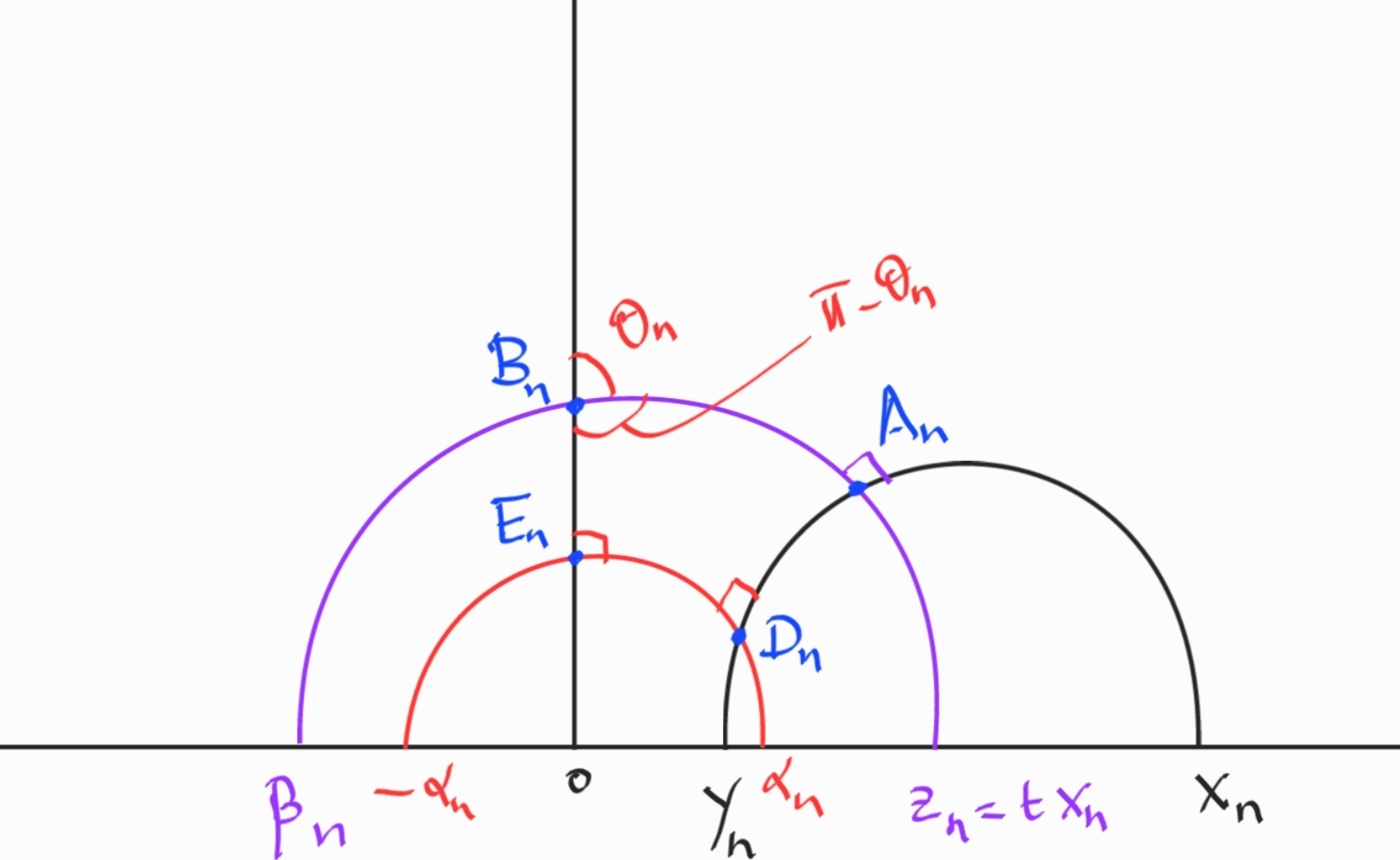}
    \caption{Polygone Hyperbolique $A_{n}B_{n}E_{n}D_{n}$}
    \label{TRI1}
\end{figure}

Comme $(Y_n)_{n\geq 1}$ converge vers $-\ds\frac{d}{c}$ alors pour $n$ assez grand $Y_n>0$ et donc $0<Y_n<X_n<\infty$. Ainsi les g\'eod\'esiques $(Y_nX_n)$ et $(0\infty)$ sont disjoints et en nottant $(-\alpha_n\alpha_n)$ la g\'eod\'esique orthogonale commune nous avons $[-\alpha_n, \alpha_n,Y_n,X_n]=-1$, ce qui donne $\alpha_n=\sqrt{Y_nX_n}$. Notons $E_n=i\alpha_n$ le point d'intersection de la g\'eod\'esique $(0\infty)$ avec la g\'eod\'esique  $(-\alpha_n\alpha_n)$ et par $D_n=t_n+is_n$ celui de la g\'eod\'esiques $(Y_nX_n)$ avec la g\'eod\'esique $(-\alpha_n\alpha_n)$, le quadrilat\`ere hyperbolique $A_nB_nE_nD_n$ a trois angles droits  et l'angle au sommet $B_n$ vaut $\pi-\theta_n$ (voir figure \ref{TRI1}).
Par les formules de trigonom\'etries hyperbolique nous avons :
 $$\sin(\pi-\theta_n)=-\sin\theta_n=\frac{\cosh(d(A_n,D_n))}{\cosh(d(B_,E_n))}.$$

 D'une part nous avons $d(B_n,E_n)=\ln\bigg(\sqrt{\ds\frac{-\beta_nZ_n}{Y_nX_n}}\bigg)=\ln\bigg(\sqrt{\ds\frac{-t\beta_n}{Y_n}}\bigg)$ qui tend vers $+\infty$ lorsque $n$ tend vers $+\infty$ et   par cons\'equent : $$\cosh(d(B_n,E_n))=\cosh\bigg(\ln\bigg(\sqrt{\ds\frac{-t\beta_n}{Y_n}}\bigg)\bigg)=\ds\frac{1}{2}\bigg(\sqrt{\ds\frac{-t\beta_n}{Y_n}}+\sqrt{\ds\frac{Y_n}{-t\beta_n}}\bigg)\approx\ds\frac{1}{2}\sqrt{\ds\frac{-t\beta_n}{Y_n}}.$$
  
  D'autre part nous avons $d(A_n, D_n))=2 Argsh\bigg(\ds\frac{|A_n-D_n|} {2\sqrt{v_ns_n}}\bigg)$. Or,
$$\frac{|A_n-D_n|^2} {4v_ns_n}=\frac{(u_n-t_n)^2+(v_n-s_n)^2}{4v_ns_n}=\frac{u_n^2+v_n^2+t_n^2+s_n^2- 2u_nt_n-2v_ns_n}{4v_ns_n}.$$ 
Comme $D_n$ appartient \`a la g\'eod\'esique $(-\alpha_n\alpha_n)$, nous avons $t_n^2+s_n^2=X_nY_n$ et :
$$\frac{|A_n-D_n|^2} {4v_ns_n}=\frac{u_n^2+v_n^2- 2u_nt_n+X_nY_n}{4v_ns_n}-\frac{1}{2}.$$ 
Mais nous avons aussi $u_n^2+v_n^2+X_nY_n=u_n(X_n+Y_n)$ car $A_n$ appartient \`a la g\'eod\'esique $(Y_nX_n)$ et donc 
$$\frac{|A_n-D_n|^2} {4v_ns_n}=\frac{u_n(X_n+Y_n)- 2u_nt_n}{4v_ns_n}-\frac{1}{2}=\frac{u_n}{4v_n}\frac{X_n+Y_n-2t_n}{s_n}-\frac{1}{2}.$$ 
En traduisant $D_n$ point d'intersection des g\'eod\'esiques $(Y_nX_n)$ et $(-\alpha_n\alpha_n)$, nous obtenons : $$t_n=\ds\frac{2X_nY_n}{X_n+Y_n}\ \mbox{et}\  s_n=\ds\frac{(X_n-Y_n)\sqrt{X_nY_n}}{X_n+Y_n}.$$ Nous d\'eduisonns :
$$\frac{|A_n-D_n|^2} {4v_ns_n}=\frac{u_n}{4v_n}\frac{X_n-Y_n}{\sqrt{X_nY_n}}-\frac{1}{2}=R_n^2\ \mbox{et}\  \cosh(d(A_n,D_n))=\cosh\bigg(2Argsh(R_n)\bigg).$$  Puisque $R_n$ tend vers $+\infty$ lorsque $n$ tend vers $+\infty$, nous avons  $$ \cosh(d(A_n,D_n))=\cosh\bigg(2Argsh(R_n)\bigg)\approx  2R_n^2\approx\frac{u_n}{2v_n}\sqrt{\frac{X_n}{Y_n}}$$ et nous d\'eduisons 
$$\sin\theta_n=-\frac{\cosh(d(A_n,D_n))}{\cosh(d(B_,E_n))}\approx-\frac{u_n}{v_n}\sqrt{\frac{X_n}{-t\beta_n}};$$ ce qui donne $\ds\lim_{n\rightarrow+\infty}\sin\theta_n=-1$ et donc $\ds\lim_{n\rightarrow+\infty}\cos\theta_n=0\neq\frac{t}{t-1}$ car $0<t<\ds\frac{1}{2}.$ Ce qui serait une contradiction.

{\it Supposons $-\ds\frac{d}{c}=0$.}

Dans ce cas, il existe une sous suite de $(Y_n){n\geq 1}$ non constante que nous noterons encore $(Y){n\geq 1}$ qui converge vers $0$ et qui est strictement positive ou strictement n\'egative.
Ainsi nous serons dans l'un des deux cas pr\'ec\'edent et le r\'esultat en d\'ecoule.
$$\hfil\Box$$

{\it D\'emonstration du th\'eor\`eme \ref{pa}}.
D 'apr\`es les lemmes \ref{lem1}, \ref{lem2} et la proposition \ref{prop1}, il existe une suite $(\gamma_n){n\geq 1}$ d'\'el\'element de $\Gamma$ telle qu'en notant $a_n, b_n, c_n$ et $d_n$ les coefficients de l'isom\'etrie $\gamma_n$, nous avons les suites $(c_n){\geq 1}$ et  $(d_n)_{n\geq 1}$ qui convergent respectivement vers $0$ et $d \in\R^{*}.$ 
Comme $\gamma_n^{-1}(\infty)=\ds\frac{-d_n}{c_n}$ et $Im(\gamma_n(i))=\ds\frac{1}{c_n^2+d_n^2}$, nous avons : $$\lim_{n\rightarrow +\infty}\gamma_n^{-1}(\infty)=\infty\ \mbox{et}\ \lim_{n\rightarrow +\infty}B_{\infty}(\gamma_n(i),i)=\lim_{n\rightarrow +\infty}\ln\bigg(c_n^2+d_n^2\bigg)=\ln(d^2) .$$ Si $d^2=1,$ le corrollaire \ref{or} assure que $h_{\R}(u)$ est r\'ecurente si non le corrollaire $\ref{tunv}$ dit que $\overline{h_{\R}(u)}$ n'est pas $h_{\R}$-minimal.

$$\hfil\Box$$

\end{document}